\documentclass[12pt]{article}
\usepackage{amscd,amssymb,epic,amsmath}

\usepackage[dvips]{graphicx}
\usepackage[mathscr]{eucal}
\title{\Large \bf On a q-analog of the Wallach-Okounkov
formula}\author {O.Bershtein, Ye.Kolisnyk, L.Vaksman}

\begin{document}
\newtheorem{proposition}{Proposition}
\newtheorem{lemma}{Lemma}
\newtheorem{corollary}{Corollary}
\newtheorem{theorem}{Theorem}

\maketitle

\centerline{Institute for Low Temperature Physics and Engineering, 47
Lenin ave. 61103, Kharkov, Ukraine.} \centerline{e-mail:
bershtein@ilt.kharkov.ua, vaksman@ilt.kharkov.ua}

\begin{abstract}
We obtain a $q$-analog of the well known result on a joint spectrum
of invariant differential operators with polynomial coefficients on a
prehomogeneous vector space of complex $n \times n$-matrices. We are
motivated by applications to the problems of harmonic analysis in the
quantum matrix ball: our main theorem can be used while proving the
Plancherel formula (to be published).


\medskip


\end{abstract}

\section{Introduction}

Our main goal is a $q$-analog of the well known Wallach-Okounkov formula
(\ref{classic_y_spectr}).

Start with recalling some well-known facts. Denote by $\mathrm{Mat}_n$ the
vector space of complex $n \times n$-matrices. The group $K=S(GL_n \times
GL_n)$ acts on $\mathrm{Mat}_n$ by
$$
(u,v)Z=uZv^{-1}, \quad (u,v) \in K, \,Z \in \mathrm{Mat}_n.
$$

This induces the $K$-action in the space $\mathbb C[\mathrm{Mat}_n]$ of
holomorphic polynomials on $\mathrm{Mat}_n$. A well-known Hua's result
claims that $\mathbb C[\mathrm{Mat}_n]$ is a direct sum of $K$-isotypic
components as follows:
$$
\mathbb C[\mathrm{Mat}_n]=\bigoplus_{\lambda \in \Lambda_n} \mathbb{C}
[\mathrm{Mat}_n]_{\bf{\lambda}}, \quad \Lambda_n=\{\lambda=
(\lambda_1,\lambda_2,\ldots,\lambda_n)\in \mathbb{Z}^n_+ | \lambda_1
\geq \lambda_2 \geq \ldots \geq \lambda_n\},
$$
where $\mathbb{C} [\mathrm{Mat}_n]_{\bf{\lambda}}$ is a simple
$K$-module with the highest weight
\begin{equation}\label{h_weight}
(\lambda_1-\lambda_2,...,\lambda_{n-1}-\lambda_n,2\lambda_n,
\lambda_{n-1}-\lambda_n,...,\lambda_1-\lambda_2).
\end{equation}

Let $Z=(z^i_j)$ be the coordinate system on $\mathrm{Mat}_n$ and
$\partial_{ij}=\frac{\partial}{\partial z^i_j}$ be the associated
differential operators on $\mathrm{Mat}_n$. Given
$$
I=\{i_1<i_2<...<i_k\} \subset \{1,...,n\}, \qquad J=\{j_1<j_2<...<j_k\} \subset
\{1,...,n\},
$$ denote by $z^I_J$ the minor of $Z=(z^i_j)$, which
line numbers are from $I$ and column numbers from $J$. Let
$\partial^I_J$ be a similar minor of $(\partial_{ij})$. Put
$$
\mathrm{y}_k=\sum_{\{I,J| \,\#(I)=\#(J)=k\}} z^I_J \cdot
\partial^I_J.
$$

$y_k|_{\mathbb{C}[\mathrm{Mat}_n]_{\mathbf \lambda}}$ is a scalar
operator, since $\mathbb{C}[\mathrm{Mat}_n]_{\mathbf \lambda}$ is a
simple $K$-module and $y_k$ is $K$-invariant. There is an explicit
formula for these scalars \cite{Wallach_operators, Sahi_to_Kostant,
Okoun}, \cite[Proposition 3.3]{KnopSahi} (the so called
Wallach-Okounkov formula):
\begin{equation}\label{classic_y_spectr}
\mathrm{y}_{k}|_{\mathbb{C}[\mathrm{Mat}_n]_{\mathbf \lambda}}\quad
=\mathfrak{s}_{\mathbf{1}^k}(\lambda_1 + n - 1,\lambda_2 + n -
2,\ldots,\lambda_{n-1}+1,\lambda_n), \qquad k=1,\ldots,n,
\end{equation}
where $\mathbf{1}^k=(\underbrace{1,...,1}_k,0,...,0)$ and the {\it
factorial} Schur polynomial $\mathfrak{s}_{\nu}$ associated to a
partition $\nu=(\nu_1,...,\nu_n)$, is defined by
$$\mathfrak{s}_\nu(x_1,x_2,\ldots,x_n) \quad=\quad
 \frac{\det\left(\prod\limits_{m=0}^{\nu_j+n-j-1}(x_i-m)
\right)_{1 \leq i,j \leq n}} {\prod\limits_{i<j} (x_i - x_j)},
$$
(see \cite{Macdonald_Schur}).

This paper presents a $q$-analog of the Wallach-Okounkov formula.

\section{The main statement}\label{Pol_Mat}

Let $q \in (0,1)$.  All algebras are assumed associative and unital
and $\mathbb C$ is the ground field.

Recall that $U_q\mathfrak{sl}_{2n}$ is a Hopf algebra with
generators $\{E_i,\:F_i,\:K_i,\:K_i^{-1}\}_{i=1}^{2n-1}$ and
relations

$$K_iK_j \,=\,K_jK_i,\quad K_iK_i^{-1}\,=\,K_i^{-1}K_i \,=\,1;$$
$$K_iE_i \,=\,q^{2}E_iK_i,\quad K_iF_i \,=\,q^{-2}F_iK_i;$$
$$K_iE_j \,=\,q^{-1}E_jK_i,\quad K_iF_j \,=\,qF_jK_i,\quad |i-j|=1;$$
$$K_iE_j \,=\,E_jK_i,\quad K_iF_j \,=\,F_jK_i,\quad |i-j| \,> \,1;$$
$$E_iF_j \,-\,F_jE_i \:=\:\delta_{ij}\frac{K_i-K_i^{-1}}{q-q^{-1}};$$
$$
E_i^2E_j\,-\,(q+q^{-1})E_iE_jE_i \,+\,E_jE_i^2 \:=\:0,\quad |i-j| \,=\,1;
$$
$$
F_i^2F_j \,-\,(q+q^{-1})F_iF_jF_i \,+\,F_jF_i^2\:=\:0, \quad |i-j|\,=\,1;
$$
$$E_iE_j-E_jE_i\,=\,F_iF_j-F_jF_i\,=\,0, \quad |i-j|\,> \,1.$$

The coproduct, the counit, and the antipode are defined as follows:
\begin{align*}
\triangle{E_j}&=E_j \otimes 1+K_j \otimes E_j,& \varepsilon(E_j)&=0,&
S(E_j)&=-K_j^{-1}E_j,\\ \triangle{F_j}&=F_j \otimes K_j^{-1}+1 \otimes F_j,&
\varepsilon(F_j)&=0,& S(F_j)&=-F_jK_j,\\ \triangle{K_j}&=K_j \otimes K_j,&
\varepsilon(K_j)&=1,& S(K_j)&=K_j^{-1},\qquad j=1,...,2n-1.
\end{align*}

Equip the Hopf algebra $U_q\mathfrak{sl}_{2n}$ with an involution $*$:
$$
(K_j^{\pm 1})^* = K_j^{\pm 1}, \quad E_j^* = \left\{
\begin{array}{rl}
   K_j F_j, & j \neq n \\
  -K_j F_j, & j = n
\end{array}
\right., \quad F_j^* = \left\{
\begin{array}{rl}
   E_j K_j^{-1}, & j \neq n \\
  -E_j K_j^{-1}, & j = n
\end{array}
\right..
$$

$U_q\mathfrak{su}_{n,n} \stackrel{\mathrm{def}}{=}
(U_q\mathfrak{sl}_{2n}, *)$ is a $*$-Hopf algebra. Denote by $U_q
\mathfrak{k} \subset U_q \mathfrak{sl}_{2n}$ the Hopf subalgebra
generated by $E_j, F_j, \, j \neq n,$ and $K_i, K_i^{-1},
i=1,...,2n-1$.

\medskip

Introduce a $*$-algebra $\mathrm{Pol}(\mathrm{Mat}_n)_q$, following
\cite{polmat}. First, denote by $\mathbb{C}[\mathrm{Mat}_n]_q$ a
well-known algebra with generators $z_a^\alpha$, $a,\alpha=1,...,n,$
and relations
\begin{align}
& z_a^\alpha z_b^\beta-qz_b^\beta z_a^\alpha=0, & a=b \quad \& \quad
\alpha<\beta,& \quad \text{or}\quad a<b \quad \& \quad \alpha=\beta,
\label{zaa1}
\\ & z_a^\alpha z_b^\beta-z_b^\beta z_a^\alpha=0,& \alpha<\beta \quad
\&\quad a>b, & \label{zaa2}
\\ & z_a^\alpha z_b^\beta-z_b^\beta z_a^\alpha-(q-q^{-1})z_a^\beta
z_b^\alpha=0,& \alpha<\beta \quad \& \quad a<b.& \label{zaa3}
\end{align}
We call $\mathbb{C}[\mathrm{Mat}_n]_q$ the algebra of holomorphic
polynomials on the quantum matrix space. It is evident that
$\mathbb{C}[\mathrm{Mat}_n]_q$ becomes the algebra of holomorphic
polynomials on complex $n \times n$-matrices under the formal
passage to the limit as $q \rightarrow 1$.

Similarly, denote by $\mathbb{C}[\overline{\mathrm{Mat}}_n]_q$ an algebra
with generators $(z_a^\alpha)^*$, $a,\alpha=1,\dots,n$ and relations
\begin{flalign}
& (z_b^\beta)^*(z_a^\alpha)^* -q(z_a^\alpha)^*(z_b^\beta)^*=0, & a=b
\quad \& \quad \alpha<\beta, & \quad \text{or} \quad a<b \quad \&
\quad \alpha=\beta, \label{zaa1*}
\\ & (z_b^\beta)^*(z_a^\alpha)^*-(z_a^\alpha)^*(z_b^\beta)^*=0, &
\alpha<\beta \quad \& \quad a>b, & \label{zaa2*}
\\ & (z_b^\beta)^*(z_a^\alpha)^*-(z_a^\alpha)^*(z_b^\beta)^*-
(q-q^{-1})(z_b^\alpha)^*(z_a^\beta)^*=0,& \alpha<\beta \quad \& \quad a<b. &
\label{zaa3*}
\end{flalign}

Denote by $\mathbb C[\mathrm{Mat}_n \oplus
\overline{\mathrm{Mat}}_n]_q$ an algebra with generators
$z_a^\alpha$, $(z_a^\alpha)^*$, $a,\alpha=1,\dots,n$, relations
(\ref{zaa1}) -- (\ref{zaa3*}), and additional relations
\begin{equation}\label{z*z}
(z_b^\beta)^*z_a^\alpha = q^2 \sum\limits_{a',b'=1}^n \sum\limits_{\alpha',
\beta' = 1}^m R(b,a,b',a') R(\beta, \alpha, \beta', \alpha')
z_{a'}^{\alpha'} \left( z_{b'}^{\beta'} \right)^* + (1-q^2) \delta_{ab}
\delta^{\alpha \beta},
\end{equation}
where $\delta_{ab}$, $\delta^{\alpha \beta}$ are Kronecker symbols,
$$
R(j,i,j',i') = \left\{
\begin{array}{cl}
  q^{-1}, &\quad i \neq j\ \&\ j=j'\ \&\ i=i', \\
  1, &\quad i=j=i'=j', \\
  -(q^{-2}-1), &\quad i=j\ \&\ i'=j'\ \&\ i'>i, \\
  0, &\quad \mbox{otherwise.}
\end{array}
\right.
$$

Finally, let $\mathrm{Pol}(\mathrm{Mat}_n)_q \stackrel{\rm def}{=}
(\mathbb C[\mathrm{Mat}_n \oplus \overline{\mathrm{Mat}}_n]_q,*)$ be
a $*$-algebra with the involution: $*: z_a^\alpha \mapsto
(z_a^\alpha)^*$. Note that our definition of
$\mathrm{Pol}(\mathrm{Mat}_n)_q$ allows us to equip it with the $U_q
\mathfrak{su}_{n,n}$-module algebra structure using Proposition
\ref{emb} and \eqref{g-act}, see the next Section.

It is very important for our goals that
$\mathrm{Pol}(\mathrm{Mat}_n)_q$ is a $q$-analog of the algebra of
differential operators with polynomial coefficients considered above.
Indeed, the latter algebra is derivable from
$\mathrm{Pol}(\mathrm{Mat}_n)_q$ via the change of generators
$z_a^\alpha \rightarrow (1-q^2)^{-1/2}z_a^\alpha$ and a subsequent
formal passage to the limit as $q \rightarrow 1$.

\bigskip
Introduce an irreducible $*$-representation of
$\mathrm{Pol}(\mathrm{Mat}_n)_q$ in a pre-Hilbert space. Denote by
$\mathcal H$ a $\mathrm{Pol}(\mathrm{Mat}_n)_q$-module with a
generator $v_0$ and defining relations
$$(z_a^{\alpha})^*v_0=0, \quad a, \alpha=1,...,n.$$

Denote by $T_F$ the representation of $\mathrm{Pol}(\mathrm{Mat}_n)_q$
which corresponds to $\mathcal H.$ Statements of the following proposition
are proved in \cite{polmat}.
\begin{proposition}\label{fock}
\begin{enumerate}
\item $\mathcal {H}=\mathbb C[\mathrm{Mat}_n]_qv_0$.
 \item $\mathcal{H}$ is
a simple $\mathrm{Pol}(\mathrm{Mat}_n)_q$-module. \item There exists a
unique sesquilinear form $(\cdot,\cdot)$ on $\mathcal H$ with the following
properties:\\ i) $(v_0,v_0)=1$; ii) $(fv,w)=(v,f^*w)$ for all $v,w \in
\mathcal H$, $f \in \mathrm{Pol}(\mathrm{Mat}_n)_q$. \item The form
$(\cdot,\cdot)$ is positive definite on $\mathcal H$. \item $T_F$ is a
faithfull representation.
\end{enumerate}
\end{proposition}

Similarly to the classical case, $\mathbb{C} [\mathrm{Mat}_n]_q =
\bigoplus_{\lambda \in \Lambda_n} \mathbb{C} [\mathrm{Mat}_n]_{q,
\bf{\lambda}},$ with $\mathbb{C} [\mathrm{Mat}_n]_{q, \bf{\lambda}}$ being
a simple $U_q \mathfrak{k}$-module with the highest weight
\eqref{h_weight}.\footnote{$\mathbb{C} [\mathrm{Mat}_n]_{q, \bf \lambda}$
is a module of type 1, see \cite{Jantzen}.} So, $\mathcal{H}$ inherits the
decomposition
\begin{equation}\label{decomp}
\mathcal{H}=\bigoplus_{\lambda \in \Lambda_n}
\mathcal{H}_{\bf{\lambda}}.
\end{equation}
Define elements $y_k \in \mathrm{Pol}(\mathrm{Mat}_n)_q$:
\begin{equation}\label{ball_y_i}
y_k = \sum\limits_{\{J'\ |\ \mathrm{card}(J') = k\}}\;
\sum\limits_{\{J''\ |\ \mathrm{card}(J'') = k\}} z_{\quad
J''}^{\wedge k\; J'} \left(z_{\quad J''}^{\wedge k\; J'}\right)^*,
\qquad k=1,\ldots,n,
\end{equation}

with $q$-minors $$ z_{\quad I}^{\wedge k\,
J}\stackrel{\mathrm{def}}{=}\sum_{s \in
S_k}(-q)^{l(s)}z_{i_1}^{j_{s(1)}} z_{i_2}^{j_{s(2)}}\cdots
z_{i_k}^{j_{s(k)}},
$$
\centerline{$I=\{(i_1,i_2,\dots,i_k)|1 \le i_1<i_2<\dots<i_k \le
n\}$,} \centerline{$J=\{(j_1,j_2,\dots,j_k)|1 \le j_1<j_2<\dots<j_k
\le n\}$.}

\begin{proposition}
$y_1,...,y_n$ are $U_q \mathfrak{k}$-invariant and $y_iy_j=y_jy_i$
for $i,j=1,...,n$.
\end{proposition}
{\bf Proof.} $U_q \mathfrak{k}$-invariance follows from explicit
calculations, while commutativity is deduced from the faithfullness of
$T_F$ and the simplicity of summands in \eqref{decomp}.\footnote{This proof
belongs to D. Shklyarov.} \hfill $\square$

As in the classical case, $T_F(y_k)|_{\mathcal{H}_\lambda}$ are scalar
operators for all $k$ and $\lambda$. Our goal is to obtain an explicit
formula for scalars $T_F(y_k)|_{\mathcal{H}_\lambda}$ (see Theorem
\ref{spectr_y_k}).

\bigskip
Recall the notation for $q$-factorial Schur polynomials \cite{KnopSahi}
$$
\mathfrak{s}_\nu(x_1,x_2,\ldots,x_n;q)=
\frac{\det\left(\prod\limits_{m=0}^{\nu_j+n-j-1}(x_i-q^m)) \right)_{1 \leq
i,j \leq n}} {\prod\limits_{i<j} (x_i - x_j)}.
$$

\begin{theorem}\label{spectr_y_k}
For all $k=1,2,\ldots,n$
\begin{equation}\label{like_Capelli_new}
T_F(y_k)|_{\mathcal{H}_\lambda}= \mathrm{const}\,\mathfrak{s}_{\mathbf{1}^k}
  (q^{2(\lambda_1 + n- 1)}, q^{2(\lambda_2 + n - 2)}, \ldots,
 q^{2(\lambda_{n-1} + 1)},  q^{2\lambda_n};q^2),
\end{equation}
with
\begin{equation}\label{const_Capelli}
\mathrm{const}=(-1)^k q^{-k(k-1)-2k(n-k)}. \end{equation}
\end{theorem}

Describe the proof briefly. Start with two auxiliary propositions.
\begin{proposition}\label{mu_lambda_new}
There exists a map $\mu: \Lambda_n \to \Lambda_n$, such that for all
$\lambda \in \Lambda_n$, $k=1,2,\ldots,n$
\begin{equation*}
T_F(y_k)|_{\mathcal{H}_\lambda}= \operatorname{const}\,
\mathfrak{s}_{\mathbf{1}^k}(q^{2(\mu(\lambda)_1+n-1)},q^{2(\mu(\lambda)_2+n-2)},
\ldots, q^{2(\mu(\lambda)_{n-1}+1)}, q^{2\mu(\lambda)_n};q^2)
\end{equation*}
with $\mathrm{const}$ as in \eqref{const_Capelli}.
\end{proposition}

\begin{proposition}\label{spec_y1}
For all $\lambda \in \Lambda_n$
$$
T_F(y_1)|_{\mathcal{H}_\lambda}=-q^{-2(n-1)}
\mathfrak{s}_{\mathbf{1}^1}(q^{2(\lambda_1+n-1)},q^{2(\lambda_2+n-2)},
\ldots,q^{2(\lambda_{n-1}+1)}, q^{2\lambda_n};q^2).
$$
\end{proposition}

Theorem 1 follows from Propositions \ref{mu_lambda_new} and
\ref{spec_y1}. Indeed, it is enough to prove it for transcendental
$q$. Prove that $\mu =\mathrm{id}$.  One can observe that the
equation
$\mathfrak{s}_{\mathbf{1}^1}(q^{2(\mu(\lambda)_1+n-1)},\ldots,
q^{2\mu(\lambda)_n};q^2)= \mathfrak{s}_{\mathbf{1}^1}
(q^{2(\lambda_1+n-1)},\ldots,q^{2\lambda_n};q^2)$ leads to
$\mu=\mathrm{id}$ for transcendental $q$.

Now we have to prove Propositions \ref{mu_lambda_new} and
\ref{spec_y1}. The proof of proposition \ref{mu_lambda_new} splits
into several steps (Subsections \ref{1} -- \ref{3}). Specifically,
Subsection \ref{1} contains a result which allows one to reduce
finding the joint spectrum of $y_k$ in the $C^*$-enveloping algebra
of $\mathrm{Pol}(\mathrm{Mat}_n)_q$ to describing the joint spectrum
of some special elements $x_k$. In Subsections \ref{2} and \ref{3} we
deduce the joint spectrum of $x_k$ from a Stokman-Dijkhuizen result
concerning an invariant integral on the quantum Grassmanian.

The proof of proposition \ref{spec_y1} is relatively simple, see
Subsection \ref{4}. Note that Proposition \ref{spec_y1} is a special
case of Theorem \ref{spectr_y_k}.

\section{Sketch of the proof}\label{result}
\subsection{From $\{y_k\}$ to $\{x_k\}$}\label{1}

Denote by $\mathbb{C}[SL_{2n}]_q$ the well-known Hopf algebra with
generators $\{t_{ij}\}_{i,j=1,\ldots,2n}$ and relations
\begin{flalign*}
& t_{\alpha a}t_{\beta b}-qt_{\beta b}t_{\alpha a}=0, & a=b \quad \& \quad
\alpha<\beta,& \quad \text{or}\quad a<b \quad \& \quad \alpha=\beta,
\\ & t_{\alpha a}t_{\beta b}-t_{\beta b}t_{\alpha a}=0,& \alpha<\beta \quad
\&\quad a>b,
\\ & t_{\alpha a}t_{\beta b}-t_{\beta b}t_{\alpha a}-(q-q^{-1})t_{\beta a}
t_{\alpha b}=0,& \alpha<\beta \quad \& \quad a<b,
\\ & \det \nolimits_q \mathbf{t}=1.
\end{flalign*}
Here $\det_q \mathbf{t}$ is a $q$-determinant of the matrix
$\mathbf{t}=(t_{ij})_{i,j=1,\ldots,2n}$:
\begin{equation*}\label{qdet}
\det \nolimits_q\mathbf{t}\stackrel{\mathrm{def}}{=}\sum_{s \in
S_{2n}}(-q)^{l(s)}t_{1\,s(1)}t_{2\,s(2)}\ldots t_{2n\,s(2n)},
\end{equation*}
with $l(s)=\mathrm{card}\{(i,j)|\;i<j \; \& \; s(i)>s(j) \}$. The
comultiplication $\Delta$, the counit $\varepsilon$, and the antipode
$S$ are defined as follows:
$$
\Delta(t_{ij})=\sum_kt_{ik}\otimes t_{kj},\qquad
\varepsilon(t_{ij})=\delta_{ij},\qquad S(t_{ij})=(-q)^{i-j}\det \nolimits_q
\mathbf{t}_{ji}
$$
with $\mathbf{t}_{ji}$ being the matrix derived from $\mathbf{t}$ by
discarding its $j$-th row and $i$-th column.

Equip $\mathbb C[SL_{2n}]_q$ with the standard $U_q
\mathfrak{sl}_{2n}$-module algebra structure as follows (see
\cite{polmat}): for $k=1,..,2n-1$
\begin{align}\label{g-act}
E_k\cdot t_{ij}= q^{-1/2} & \begin{cases} t_{i\,j-1}, & k=j-1,
\\ 0, & \text{otherwise},
\end{cases} \qquad
F_k \cdot t_{ij}= q^{1/2}
\begin{cases}
t_{i\,j+1}, & k=j,
\\ 0, & \text{otherwise},
\end{cases}
\\ & K_k \cdot t_{ij}=
\begin{cases}
qt_{ij}, & k=j,
\\q^{-1}t_{ij}, & k=j-1,
\\ t_{ij}, & \text{otherwise}.
\end{cases}
\end{align}
Let $\mathbb{C}[X]_q \stackrel{\mathrm{def}} {=}
(\mathbb{C}[SL_{2n}]_q,*)$ be a $*$-algebra with the involution $*$
given by
\begin{equation}\label{*_noncom}
t_{ij}^*= \mathrm{sign} (i-n-1/2)(n-j+1/2)(-q)^{j-i}\det \nolimits_q
\mathbf{t}_{ij}.
\end{equation}
It is a $U_q \mathfrak{su}_{n,n}$-module $*$-algebra.
 Recall a standard notation for $q$-minors of $\mathbf{t}$:
$$
t_{IJ}^{\wedge k}\stackrel{\mathrm{def}}{=}\sum_{s \in
S_k}(-q)^{l(s)}t_{i_1j_{s(1)}}\cdot t_{i_2j_{s(2)}}\cdots
t_{i_kj_{s(k)}},
$$
with $I=\{1 \le i_1<i_2<\dots<i_k \le 2n\}$, $J=\{1 \le j_1<j_2<\dots<j_k
\le 2n\}$. Introduce the elements
$$
t=t_{\{1,2,\dots,n \}\{n+1,n+2,\dots,2n \}}^{\wedge n},\qquad x=tt^*.
$$
Note, that $t,$ $t^*,$ and $x$ quasi-commute with all generators
$t_{ij}$ of $\mathbb{C}[SL_{2n}]_q$, and $\mathbb{C}[X]_q$ is an
integral domain \cite{BrGood}. Let $\mathbb{C}[X]_{q,x}$ be the
localization of $\mathbb{C}[X]_q$ with respect to the multiplicative
set $x^{\mathbb Z_+}$. The following statements are proved, for
instance, in \cite{polmat}.

\begin{proposition}\label{local}
There exists a unique extension of the $U_q
\mathfrak{su}_{n,n}$-module $*$-algebra structure from
$\mathbb{C}[X]_q$ onto $\mathbb{C}[X]_{q,x}$.
\end{proposition}

\begin{proposition}\label{emb}
The map
\begin{equation}
i:z_a^\alpha \mapsto t^{-1}t_{\{1,2,\dots,n \}J_{a \alpha}}^{\wedge n},
\end{equation}
with $J_{a \alpha}=\{n+1,n+2,\dots,2n \}\setminus \{2n+1-\alpha
\}\cup \{a \}$, admits a unique extension up to an embedding of $U_q
\mathfrak{su}_{n,n}$-module $*$-algebras
$i:\mathrm{Pol}(\mathrm{Mat}_n)_q \hookrightarrow
\mathbb{C}[X]_{q,x}$.
\end{proposition}

The last proposition allows to identify
$\mathrm{Pol}(\mathrm{Mat}_n)_q$ with its image in
$\mathbb{C}[X]_{q,x}$. It can be proved that
$$y_k=(-1)^k(tt^*)^{-1}\sum_{\substack{K \subset
\{1,...,2n\}, \;\mathrm{card}(K)=n,\\
  \mathrm{card}(K \cap \{1,2,\ldots,n\})=k}}
(-q)^{l(K,K^c)} \; t_{\{1,2,\ldots,n \}\,K}^{\wedge n}\;
t_{\{n+1,n+2,\ldots,2n \}\, K^c}^{\wedge n},$$ where
$K^c=\{1,...,2n\} \setminus K$ and $l(K,K^c)=\mathrm{card}((i,j)| i
\in K, j \in K^c, i>j)$. Consider elements of $\mathbb{C}[X]_q$:
$$x_k\,=\,q^{k(k-1)} \sum_{\substack{I \subset \{1,2,\ldots,n\},\, J\subset
\{n+1,n+2,\ldots,2n\} \\ \mathrm{card} (I)= \mathrm{card} (J)=k}} q^{-2
\sum\limits_{m=1}^k (n-i_m)} (-q)^{\sum\limits_{m=1}^k (j_m-i_m-n)}\,
t_{I\,J}^{\wedge k} \, t_{I^c\,J^c}^{\wedge (2n-k)}.$$

  We substitute the problem of computing the joint spectrum of $y_k$ with a problem
of computing the joint spectrum of $x_k$ via the next proposition.
Put $\begin{pmatrix} a \\
b \end{pmatrix}_{q} = \frac{(q;q)_a}{(q;q)_b (q;q)_{a-b}}$.

\begin{proposition}\label{x_y_link} For all $k=1,2,\ldots,n$
we have
\begin{equation*}
x_k = \frac{\sum\limits_{m=0}^{n-k}(-1)^m\left(\begin{array}{c} n-m \\
k \end{array}\right)_{q^{-2}} y_m} {\sum\limits_{m=0}^n (-1)^my_m}
\end{equation*} with  $y_0=1$.
\end{proposition}
The proof of this statement can be managed by explicit computations
in $\mathbb C[SL_{2n}]_q$ and is omitted.

\subsection{From $*$ to $\star$}\label{2}

In this subsection we suppose that $q>0$ and $q \neq 1$ instead of $q \in
(0,1)$. Introduce a $*$-Hopf algebra $U_q\mathfrak{su}_{2n}
\stackrel{\mathrm{def}}{=} (U_q\mathfrak{sl}_{2n}, \star)$ with an
involution $\star$:
$$
(K_j^{\pm 1})^\star = K_j^{\pm 1}, \quad E_j^\star =
 K_j F_j,
\quad F_j^\star = E_j K_j^{-1}.
$$

Let
$\mathbb{C}[SU_{2n}]_q\stackrel{\mathrm{def}}{=}(\mathbb{C}[SL_{2n}]_q,\star)$
be a $*$-Hopf algebra with an involution $\star$ given by
\begin{equation}\label{*_com}
t_{ij}^\star= (-q)^{j-i}\det \nolimits_q \mathbf{t}_{ij}.
\end{equation}
It is well-known that $\mathbb{C}[SU_{2n}]_q$ is a $U_q
\mathfrak{su}_{2n}$-module algebra.

Similarly to Proposition \ref{local}, the structure of $U_q
\mathfrak{su}_{2n}$-module algebra extends up to the localization
$\mathbb{C}[SU_{2n}]_{q,x}$ of the algebra $\mathbb{C}[SU_{2n}]_q$ with
respect to the multiplicative set $x^{\mathbb Z_+}$. Note that
\begin{align*}
&x_k\,=\,q^{k(k-1)} \sum_{\substack{I \subset \{1,2,\ldots,n\},\,
J\subset \{n+1,n+2,\ldots,2n\}\\
\mathrm{card}(I)=\mathrm{card}(J)=k}} q^{-2 \sum\limits_{m=1}^k
(n-i_m)} \, t_{I\,J}^{\wedge k} \, (t_{I\,J}^{\wedge k})^*
\\ &=\,q^{k(k-1)} \sum_{\substack{I \subset \{1,2,\ldots,n\},\, J\subset
\{n+1,n+2,\ldots,2n\}\\
\mathrm{card}(I)=\mathrm{card}(J)=k}} q^{-2 \sum\limits_{m=1}^k
(n-i_m)} \, t_{I\,J}^{\wedge k} \, (t_{I\,J}^{\wedge k})^\star.
\end{align*}

Up to the end of this subsection our aim is obtaining the equality
\eqref{T_T}. Roughly speaking, we describe an interplay between the joint
spectrum of $x_k$ in a $*$-representation of $\mathbb{C}[SU_{2n}]_q$ and in
a $*$-representation of $\mathbb{C}[X]_q$.

Let $\operatorname{Pol}(U)_q \stackrel{\mathrm{def}}{=}
(\mathbb{C}[\mathrm{Mat}_n \oplus
\overline{\mathrm{Mat}}_n]_q,\star)$ be a $*$-algebra with an
involution $\star$ (see Section \ref{Pol_Mat}):
$$  f^\star =(-1)^{\rm{deg}\,f}\, f^*, \qquad f \in
\mathbb{C}[\mathrm{Mat}_n]_q,$$ where $\deg f$ is defined in a natural way.
One can verify that $\operatorname{Pol}(U)_q$ is a
$U_q\mathfrak{su}_{2n}$-module algebra.\footnote{ $\operatorname{Pol}(U)_q$
is a $q$-analog of the space of regular functions on the big cell of
Grassmanian $Gr_n(\mathbb C^{2n})$ considered as a real algebraic variety.}
 Introduce an irreducible $*$-representation of $\operatorname{Pol}(U)_q$
in a pre-Hilbert space. Let $\mathscr{H}$ be a
$\operatorname{Pol}(U)_q$-module with a single generator $v^0$ and
relations
$$ z_a^\alpha\,v^0=0,\qquad a,\alpha=1,2,\ldots, n.$$ Denote by
$\mathscr{T}_F$ the corresponding representation of
$\operatorname{Pol}(U)_q$. Similarly to Proposition \ref{fock}, we
have
\begin{proposition}
\begin{enumerate}
\item $\mathscr {H}=\mathbb C[\overline{\mathrm{Mat}}_n]_qv^0$.
 \item $\mathscr{H}$ is
a simple $\mathrm{Pol}(U)_q$-module. \item There exists a unique
sesquilinear form $(\cdot,\cdot)$ on $\mathscr H$ with the following
properties:\\ i) $(v^0,v^0)=1$; ii) $(fv,w)=(v,f^\star w)$ for all
$v,w \in \mathscr H$, $f \in \mathrm{Pol}(U)_q$; iii) the form
$(\cdot,\cdot)$ is positive definite on $\mathscr H$.
\item $\mathscr{T}_F$ is a faithfull representation.
\end{enumerate}
\end{proposition}

Equip $\mathbb{C}[\mathrm{Mat}_n]_q$ and
$\mathbb{C}[\overline{\mathrm{Mat}}_n]_q$ with the gradings:
\begin{align*}
\mathbb{C}[\mathrm{Mat}_n]_q=\bigoplus_{j=0}^\infty
\mathbb{C}[\mathrm{Mat}_n]_{q,j}, \qquad
\mathbb{C}[\mathrm{Mat}_n]_{q,j}=\{f \in
\mathbb{C}[\mathrm{Mat}_n]_q| \deg f=j\},
\\ \mathbb{C}[\overline{\mathrm{Mat}}_n]_q=\bigoplus_{j=0}^\infty
\mathbb{C}[\overline{\mathrm{Mat}}_n]_{q,-j}, \qquad
\mathbb{C}[\overline{\mathrm{Mat}}_n]_{q,-j}=\{f \in
\mathbb{C}[\overline{\mathrm{Mat}}_n]_q| \deg f=j\}.
\end{align*}

Then vector spaces $\mathcal{H}$ and $\mathscr{H}$ inherits the
corresponding gradings:
$$\mathcal{H}=\bigoplus \limits_{j=0}^\infty\mathcal{H}_j,\quad
\mathcal{H}_j= \mathbb{C}[\mathrm{Mat}_n]_{q,j}v_0 \qquad \text{and} \qquad
\mathscr{H}=\bigoplus \limits_{j=0}^\infty\mathscr{H}_j,\quad
\mathscr{H}_j= \mathbb{C}[\overline{\mathrm{Mat}}_n]_{q,-j}v^0.
$$

Consider an antilinear algebra antiautomorphism $$ \sigma:
\mathbb{C}[SL_{2n}]_q \rightarrow \mathbb{C}[SL_{2n}]_{q^{-1}},
\qquad \sigma: t_{ij} \mapsto t_{ij}.
$$
\begin{lemma}\label{x-x-const} For all $k=1,2,\ldots,n$
$$\sigma(x_k(q))\,=\,q^{2k^2}\, x_k(q^{-1}).$$
\end{lemma}
The proof reduces to explicit computations. $\sigma$ extends up to an
antihomomorphism $\mathbb{C}[SL_{2n}]_{q,x(q)} \rightarrow
\mathbb{C}[SL_{2n}]_{q^{-1},x(q^{-1})}$. Thus there exists a
$*$-representation of $\mathrm{Pol}(U)_q$
$$
\widetilde{T}_F(f)=\mathscr{T}_F(\sigma(f))^*, \qquad f \in
\mathbb{C}[\mathrm{Mat}_n \oplus \overline{\mathrm{Mat}}_n]_q
$$
in the space $\mathscr{H}^*\stackrel{\mathrm{def}}{=} \oplus_{j=0}^\infty
\mathscr{H}_j^*$.

\begin{lemma}\label{Fock_again}
$\widetilde{T}_F$ and $T_F$ are equivalent representations of
$\mathbb{C}[\mathrm{Mat}_n \oplus \overline{\mathrm{Mat}}_n]_q$.
\end{lemma}

{\bf Proof.} Each linear map $i: \mathcal{H}_0 \rightarrow
\mathscr{H}_0$ extends to a $\mathbb{C}[\mathrm{Mat}_n \oplus
\overline{\mathrm{Mat}}_n]_q$-morphism
$$
i: \mathcal{H} \rightarrow \mathscr{H}^*, \quad i(\mathcal{H}_j)=\mathscr{H}^*_j.
$$
Since $\mathrm{dim}\mathcal{H}_j=\mathrm{dim}\mathscr{H}^*_j$ and
$\mathcal{H}$ is simple, the proof is done. \hfill $\square$

The representation $T_F$ of $\mathbb{C}[\mathrm{Mat}_n \oplus
\overline{\mathrm{Mat}}_n]_q$ can be naturally considered as the
$*$-representation of the isomorphic subalgebra in
$\mathbb{C}[X]_{q,x}$ (see Prop. \ref{emb}). Similarly, the
representation $\mathscr{T}_F$ of $\mathbb{C}[\mathrm{Mat}_n \oplus
\overline{\mathrm{Mat}}_n]_q$ can be naturally considered as the
$*$-representation of the isomorphic subalgebra in
$\mathbb{C}[SU_{2n}]_{q,x}$.

By Lemmas \ref{x-x-const}, \ref{Fock_again}, we have
\begin{equation}\label{T_T}
 T_F(x_k(q))|_{\mathcal{H}_\lambda} = q^{2k^2}
\left(\mathscr{T}_F(x_k(q))|_{\mathscr{H}_\lambda} \right)_{q
\mapsto q^{-1}}, \qquad \lambda \in \Lambda_n.
\end{equation}

\subsection{Stokman-Dijkhuizen results}\label{3}

The next step in the proof of Theorem 1 is obtaining a formula for the
joint spectrum of $\{\mathscr{T}_F(x_k)\}$. We use the results of
\cite{StokmanPhd} about zonal spherical functions and invariant measure on
quantum analogs of $Gr_n(\mathbb C^{2n})$. In this paper Stokman and
Dijkhuizen introduce pairwise commuting elements $\widetilde{e^{\infty,
\infty}_1},\widetilde{e^{\infty, \infty}_2},...,\widetilde{e^{\infty,
\infty}_n}$ and find an explicit formula for the discrete part of their
joint spectrum in the $C^*$-enveloping algebra of $\mathbb C[SU_{2n}]_q$.

Denote by $e_k(z_1, \ldots, z_n)$ the elementary symmetric
polynomial of $n$ variables. In the following proposition the first
statement belongs to Stokman and Dijkhuizen \cite{StokmanPhd}, while
the second is a useful supplement.

\begin{proposition}\label{support_nu}
1. The joint spectrum of elements $\widetilde{e_k^{\infty, \infty}}$ is the
closure of the set
$$\{(e_1(\lambda),e_2(\lambda),\ldots, e_n(\lambda))\,|\, \lambda \in \Lambda_n\},$$
where $e_j(\lambda)$ stands for $e_j(q^{2(\lambda_1+n-1)},
q^{2(\lambda_2+n-2)}, \ldots,q^{2(\lambda_{n-1}+1)}, q^{2\lambda_n})$ for
short.

2. $x_k=q^{k(k+1)}\widetilde{e_k^{\infty, \infty}}$.

3. The set of joint eigenvalues of $\mathscr{T}_F(x_k)$ is
$$\{(a_1(\lambda),a_2(\lambda),\ldots, a_n(\lambda))\,|\, \lambda \in \Lambda_n\},$$ with
$a_k(\lambda)=q^{k(k+1)}e_k(q^{2(\lambda_1+n-1)}, q^{2(\lambda_2+n-2)},
\ldots,q^{2(\lambda_{n-1}+1)}, q^{2\lambda_n}),$ $k= 1,2,\ldots,n.$
\end{proposition}

Note that $a_k(\lambda)$ are polynomials of $q$. Using \eqref{T_T}, we get
\begin{corollary}\label{spectr_of_x}
There exists a map $\mu: \Lambda_n \to \Lambda_n$ such that
\begin{equation}\label{spec_x}
 T_F(x_k)|_{\mathcal{H}_\lambda} = q^{k(k-1)}
e_k(q^{-2(\mu(\lambda)_1+n-1)},q^{-2(\mu(\lambda)_2+n-2)},\ldots,
q^{-2(\mu(\lambda)_{n-1}+1)},q^{-2\mu(\lambda)_n})
\end{equation}
for all $k=1,2,\ldots,n$, $\lambda \in \Lambda_n$.
\end{corollary}

The last part of the proof of Theorem 1 includes the following identity:
\begin{proposition}\label{poly_new} In $\mathbb{C}[z_1,z_2,\ldots,z_n]$, we
have
\begin{equation*}
q^{k(k-1) - n(n-1)} e_{n-k}(z_1,\ldots,z_n) = \sum\limits_{m=0}^{n-k}
q^{-m(2n-m-1)}
\left(\begin{array}{c} n-m \\
k \end{array}\right)_{q^{-2}} s_{\mathbf{1}^m}(z_1,\ldots,z_n;q^2)
\end{equation*}
for all $k=1,2,\ldots,n$.
\end{proposition}

The proof is based on results of \cite{Knop-Capelli}. \hfill
$\square$

\subsection{Proof of Proposition \ref{spec_y1}}\label{4}

Due to the definition of $q$-factorial Schur function, Proposition
\ref{spec_y1} is equivalent to
\begin{proposition}\label{4_1}
For all $\lambda \in \Lambda_n$
\begin{equation} \label{eigenv_y_1}
 T_F(y_1)|_{\mathcal{H}_\lambda} = \sum_{j=1}^{n}
 q^{-2(j-1)}- \sum_{j=1}^{n}
 q^{2(\lambda_j-(j-1))}.
\end{equation}
\end{proposition}

The proof is based on the following easy result predicted by D.Shklyarov
\begin{lemma}\label{D_res}
For $n>1$ the map
$$
J_n: z^\alpha_a \mapsto
\begin{cases} q^{-1} z^\alpha_a, &
              \quad 1 \leq \alpha \leq n-1 \; \&\; 1\leq a \leq n-1,\\
               1, & \quad a=\alpha=n\\
               0, & \quad \text{otherwise}.
\end{cases}
$$
has a unique extension up to a $*$-algebras homomorphism $J_{n}:\;
\mathrm{Pol}(\mathrm{Mat}_n)_q \to \mathrm{Pol}
(\mathrm{Mat}_{n-1})_q.$
\end{lemma}
{\bf Proof of Proposition \ref{4_1}.} We proceed by induction in
$n$. For $n=1$ the statement is obvious.

Let
\begin{equation*}\label{u_lambda}
u_\lambda=({\rm det}_q{\mathbf z})^{\lambda_n} \cdot \prod_{j=1}^{n-1}
  \left(z_{\quad\, \{1,2,\ldots,j\}}^{\wedge j\,
\{1,2,\ldots,j\}}\right)^{\lambda_j-\lambda_{j+1}}\; v_0 \in
\mathcal{H}_\lambda.
\end{equation*}
Using Lemma \ref{D_res} and the fact that
$$
\dim (\mathbb{C}[\mathrm{Mat}_n]_{q,n+1} \otimes
\mathbb{C}[\overline{\rm Mat}_n]_{q,1}
 + \mathbb{C}[\mathrm{Mat}_n]_{q,n} \otimes 1)^{U_q\mathfrak{k} \otimes U_q\mathfrak{k}}
 \leq 2,
$$ one can easy prove that
\begin{equation}\label{coroll_1} y_1\, \mathrm{det}_q{\mathbf z}=
 q^2\, \mathrm{det}_q{\mathbf z}\;( y_1\; +\; q^{-2n} -1),
\end{equation}
so
\begin{equation}\label{det_prev}
T_F(y_1)u_\lambda= T_F((\mathrm{det}_q{\mathbf z})^{\lambda_n}) \cdot
T_F\left(q^{2\lambda_n}y_1+(1-q^{2\lambda_n})\sum_{j=1}^{n}q^{-2(j-1)}
)\right) \prod_{j=1}^{n-1}
  \left(z_{\quad\, \{1,2,\ldots,j\}}^{\wedge j\,
\{1,2,\ldots,j\}}\right)^{\lambda_j-\lambda_{j+1}}\; v_0.
\end{equation}

The subalgebra in $\mathbb{C}[{\rm Mat}_n]_q$ generated by
$z_a^\alpha$ for $\alpha,a < n$ is isomorphic to the algebra
$\mathbb{C}[{\rm Mat}_{n-1}]_q$, and is related to a subspace
$\mathcal{H}'=\mathbb{C}[{\rm Mat}_{n-1}]_q v_0$ of the pre-Hilbert
space $\mathcal{H}=\mathbb{C}[{\rm Mat}_n]_q v_0$.

It follows from \eqref{z*z} that in the proof we can restrict ourself by
the $*$-algebra $\mathrm{Pol}(\mathrm{Mat}_{n-1})_q$ and its representation
in $\mathcal{H}'$. By the induction hypothesis
$$(1-q^{2\lambda_n})\sum_{j=1}^n q^{-2(j-1)} +
q^{2\lambda_n}\left(\sum_{j=1}^{n-1}q^{-2(j-1)}-
\sum_{j=1}^{n-1}q^{2(\lambda_j -\lambda_n -(j-1))}\right)=
\sum_{j=1}^n q^{-2(j-1)} - \sum_{j=1}^{n}q^{2(\lambda_j-(j-1))}.
\eqno \square$$

\section{Acknowledgement}

The last named author thanks to G. Zhang for attracting our
attention to Wallach's paper and to J. Stokman for helpful
discussions. Also we would like to thank F. Knop for sending the
authors to Okounkov's result.

We devote the paper to our friend and colleague Dmitry Shklyarov for
the occasion of his 30-th birthday.

\end{document}